\documentclass[reqno]{amsart}

\usepackage{amsthm,amssymb,amsmath,amscd,graphicx,latexsym,enumerate}

\newcommand{\fg}{{\mathfrak g}}
\newcommand{\fh}{{\mathfrak h}}
\newcommand{\CC}{{\mathbb C}}
\newcommand{\RR}{{\mathbb R}}

\newcommand{\ZZ}{{\mathbb Z}}

\newcommand{\K}{{\mathcal K}}

\newcommand{\fk}{{\mathfrak k}}
\newcommand{\fn}{{\mathfrak n}}

\newcommand{\M}{{\mathcal M}}
\newcommand{\MM}{{\mathcal M}}

\renewcommand{\P}{{\mathcal P}}
\renewcommand{\SS}{{\mathbb S}}
\theoremstyle{definition}
\newtheorem*{definition*}{Definition}
\newtheorem{theorem}{Theorem}[section]
\newtheorem{definition}{Definition}[section]
\newtheorem*{remark*}{Remark}
\newtheorem*{remarks*}{Remarks}
\newtheorem*{theorem*}{Theorem}
\newtheorem{example}{Example}
\newtheorem{example*}{Example}
\newtheorem{lemma}{Lemma}[section]
\newtheorem{corollary}{Corollary}[section]
\newtheorem{proposition}{Proposition}[section]
\newtheorem*{acycax}{Acyclicity axiom}
\newtheorem{conjecture}{Conjecture}

\numberwithin{equation}{section}

\begin{document}

\title{Morse theory on graphs}
\author[V. Guillemin]{V. Guillemin\footnotemark {*}}
\thanks{* Supported by NSF grant DMS 890771.}
\address{Department of Mathematics, MIT, Cambridge, MA 02139}
\email{vwg@@math.mit.edu}
\author[C. Zara]{C. Zara\footnotemark {**} }
\thanks{** Supported by Alfred P. Sloan Doctoral Dissertation 
Fellowship grant DD 766.}
\address{Department of Mathematics, MIT, Cambridge, MA 02139}
\email{czara@math.mit.edu}

\begin{abstract}
Let $\Gamma$ be a finite $d$-valent graph and $G$ an $n$-dimensional 
torus. An ``action'' of $G$ on $\Gamma$ is defined by a map which assigns 
to each oriented edge, $e$, of $\Gamma$, a one-dimensional representation of 
$G$ (or, alternatively, a weight, $\alpha_e$, in the weight lattice of $G$. 
For the assignment, $e \to \alpha_e$, to be a schematic description of a 
``$G$-action'', these weights have to satisfy certain compatibility 
conditions: the GKM axioms). As in \cite{GKM} we attach to 
$(\Gamma, \alpha)$ an equivariant cohomology ring, 
$H_G(\Gamma)=H(\Gamma,\alpha)$. By definition this ring contains the 
equivariant cohomology ring of a point, $\SS(\fg^*) = H_G(pt)$, as a 
subring, and in this paper we will use graphical versions of standard 
Morse theoretical techniques to analyze the structure of $H_G(\Gamma)$ 
as an $\SS(\fg^*)$-module.
\end{abstract}

\maketitle

\begin{center}
\section{Introduction}
\end{center}

\smallskip

As Richard Stanley so astutely observes in \cite{St},  
``The number of systems of terminology presently used in graph theory 
is equal, to a close approximation, to the number of graph theorists.''
Our terminological conventions in this paper will be the following: Given
a finite $d$-valent graph, $\Gamma$, we will denote by $V_{\Gamma}$ the 
vertices of $\Gamma$ and by $E_{\Gamma}$ the oriented edges of $\Gamma$. 
For each oriented edge, $e$, we will denote by $\bar{e}$ the edge 
obtained by reversing the orientation of $e$ and by $i(e)$ and $t(e)$ the 
initial and terminal vertices of $e$. (Thus ``$d$-valent'' means that for 
every vertex, $p$, there are exactly $d$ edges with $i(e)=p$.)

Let $G$ be a torus of dimension $n>1$, and let $\varrho$ be a function which 
associates to every oriented edge, $e$, of $\Gamma$, 
a one-dimensional representation, $\varrho_e$, of $G$, 
with character $\chi_e : G  \to S^1$, 
and $\tau$ a function which associates to 
every vertex, $p$, of $\Gamma$, a $d$-dimensional representation, $\tau_p$, 
of $G$. We will say that $\varrho$ and $\tau$ define \emph{an action of} $G$ 
\emph{on} $\Gamma$ if they satisfy the three axioms below:

\noindent Axiom A1: $$ \tau_p \simeq \bigoplus_{i(e)=p} \varrho_e.$$
Axiom A2: $$\varrho_{\bar{e}} \simeq \varrho_e^*.$$
Axiom A3: 

Let $G_e$ be the kernel of $\chi_e : G \to S^1$ and let 
$p=i(e)$ and $q=t(e)$. Then 
$$\tau_p |_{G_e} \simeq \tau_q |_{G_e}.$$

The data $(\Gamma, \rho, \tau)$ can be viewed as the schematic description of 
a genuine $G$-action, namely an action of $G$ on a $d$-dimensional complex 
manifold, $M$, having the following properties:

\begin{enumerate}[{P}1]
\item :  The fixed point set of $G$ is $V_{\Gamma}$, and, for every 
$p \in V_{\Gamma}$, $\tau_p$ is the isotropy representation of $G$ on the 
tangent space to $M$ at $p$.
\item : Each edge, $e$, of $G$, corresponds to a $G$-invariant imbedded 
projective line, $\gamma : \CC P^1 \to M_e \subset M$. The points 
$p=i(e)=\gamma([1:0])$ and $q=t(e)=\gamma([0:1])$ are the two fixed points 
for the action of $G$ on $M_e$, and $\varrho_e$ is the isotropy representation 
of $G$ on the tangent space to $M_e$ at $p$.
\end{enumerate}

These properties are exhibited by a class of manifolds, called 
\emph{GKM manifolds}, which were first introduced in \cite{GKM} 
by Goresky, Kottwitz and MacPherson and have subsequently been 
studied in a number of recent 
papers (\cite{GZ1}, \cite{GZ2}, \cite{KR}, \cite{LLY}, \cite{TW}) by 
ourselves and by others. By definition, a complex $G$-manifold, $M$, is 
\emph{GKM} if, for every fixed point, $p$, of $G$, the weights, $\alpha_i$, 
$i=1,..,d$, of the isotropy representation of $G$ on the tangent space to $M$ 
at $p$ are pairwise linearly independent: \emph{i.e.} $\alpha_i$ and 
$\alpha_j$ are linearly independent if $i \neq j$. For such a manifold we 
define its \emph{one-skeleton} to be the set of points, $p \in M$, for which 
the stabilizer group of $p$ is of dimension at least $n-1$. What Goresky, 
Kottwitz and MacPherson observe is that this one-skeleton consists of closed 
submanifolds of $M$ on which $G$ acts in a fixed point free fashion and 
imbedded $\CC P^1$'s satisfying P1-P2. 

In \cite{GZ2} we proved a converse to this result. Let 
$(\Gamma, \varrho, \tau)$ be a graph-theoretical $G$-space in the sense 
described above. For $e \in E_{\Gamma}$ let $\alpha_e$ be the weight of 
the representation, $\varrho_e$, and assume that for every vertex, $p$, of 
$\Gamma$, the weights $\alpha_e$, $i(e)=p$, of the representation, $\tau_p$, 
are pairwise linearly independent. Then there exists a complex $G$-manifold, 
$M_{\Gamma}$, of GKM type, whose one-skeleton has the properties P1-P2. 
(In fact, $M_{\Gamma}$ is canonically constructed from the data 
$(\Gamma, \rho, \tau)$ and is basically just a tubular neighborhood of 
this one-skeleton.)

From Axiom A1, it is clear that the function, $\rho$, determines the function, 
$\tau$; and $\rho$, in turn, is determined by the function 
$$\alpha : E_{\Gamma} \to \fg^*, \qquad e \to \alpha_e.$$
In the geometric realization of $\Gamma$ which we described above this 
function tells us how each of the two-spheres in the one-skeleton is 
rotated about its axis of symmetry; so we will call this function the 
\emph{axial function} of the $G$-action on $\Gamma$. The axioms, A1-A3, can 
easily be translated into axioms on this axial function (as we will do in 
Section \ref{sec:2}).

By realizing $(\Gamma, \rho, \tau)$ as a $G$-manifold, $M_{\Gamma}$, we can 
attach to it a cohomology ring: the equivariant cohomology ring, 
$H_G(M_{\Gamma})$. The constant map, $M_{\Gamma} \to  pt$, makes this 
into a module over the ring, $H_G(pt)= \SS(\fg^*)$; and we will denote by 
$H_G(\Gamma)$ the torsion-free part of this module. The questions we 
will be concerned with in this paper are:

\begin{enumerate}
\item Is this module a free $\SS(\fg^*)$-module ?
\item If so, how many generators does it have in each dimension ?
\item Is there a nice combinatorial description of these generators ?
\end{enumerate}

We will now formulate these questions more precisely. We will henceforth 
abandon any pretense of being concerned with $G$-actions on 
manifolds and deal only with $G$-actions on 
graphs. Fortunately, there is a very beautiful description in \cite{GKM} 
of $H_G(\Gamma)$ in terms of $\Gamma$ and $\alpha$ which allows us 
to do this.

Fix a vector, $\xi \in \fg$, with the property that 
\begin{equation}
\label{eq:1.2}
\alpha_e(\xi) \neq 0
\end{equation}
for all $e \in E_{\Gamma}$. We will call such a vector \emph{polarizing}, and 
denote by $\P$ the set of polarizing vectors. By axiom A2
\begin{equation}
\label{eq:1.1new}
\alpha_{\bar{e}}(\xi) = - \alpha_e(\xi),
\end{equation}
so we can orient $\Gamma$ by assigning to each unoriented edge, the 
orientation which makes $\alpha_e(\xi)$ positive. We will denote this 
orientation by $o_{\xi}$. To be able to do Morse theory on $\Gamma$ we 
make the following essential assumption.

\begin{acycax}
For some $\xi$ satisfying \eqref{eq:1.2}, the orientation, $o_{\xi}$, 
of $\Gamma$ that we have just described has no oriented cycles.
\end{acycax}

We will say that $(\Gamma, \alpha)$ is \emph{$\xi$-acyclic} if this axiom 
is satisfied.

Let $V=V_{\Gamma}$. The acyclicity assumption 
implies that one can find a function $\phi : V \to \RR$ with the 
property that, for every oriented edge, $e$, of $\Gamma$, 
with initial vertex $p$ and terminal vertex $q$,
$\phi(p) < \phi(q)$. We will say that such a function 
is $\xi$-\emph{compatible} and call any function with this 
property a \emph{Morse function}. 
A canonical choice of such a $\phi$ is the following: For every vertex, 
$p$, consider all oriented paths in ($\Gamma$, $o_{\xi}$) 
with terminal  point $p$ and 
let $\phi(p)$ be the length of the longest of these paths. It is clear that 
for this function to be unambiguously defined, $\Gamma$ has to satisfy the 
acyclicity assumption.  By a small perturbation if necessary, we may 
assume, without loss of generality, that $\phi$ is injective.

Following \cite{GKM} 
we define the \emph{cohomology ring}, $H_G(\Gamma)$, of 
$(\Gamma, \alpha)$, to be the subring of the graded ring 
\begin{equation}
\label{eq:1.4}
\mbox{Maps}(V, \SS(\fg^*)) = \bigoplus_{k \geq 0}
\mbox{Maps}(V, \SS^k(\fg^*))
\end{equation}
consisting of all maps, $f : V \to \SS(\fg^*)$, which satisfy
\begin{equation}
\label{eq:1.5}
f_p \equiv f_q \pmod{\alpha_e}
\end{equation}
for every edge, $e$, of $\Gamma$, $f_p$ and $f_q$ being the values 
of $f$ at the endpoints, $p$ and $q$, of $e$. 
Since $f_p$ and $f_q$ are elements of $\SS(\fg^*)$, they can be 
thought of as polynomial functions on $\fg$ and \eqref{eq:1.5} asserts that 
the restriction of $f_p$ to the hyperplane $\alpha_e=0$ is equal to the 
restriction of $f_q$. More precisely, if $\fg_e^* = \fg^*/\RR\alpha_e$, 
then the projection $\fg^* \to \fg_e^*$ induces an epimorphism of 
rings
$$r_e : \SS(\fg^*) \to \SS(\fg_e^*)$$
and condition \eqref{eq:1.5} can be expressed as 
\begin{equation}
\label{eq:rhoe}
r_e(f_{i(e)}) = r_e(f_{t(e)}).
\end{equation}

Since the constant maps of $V$ into $\SS(\fg^*)$ satisfy 
\eqref{eq:1.5},
$\SS(\fg^*)$ is a subring of $H_G(\Gamma)$; so, in particular, 
$H_G(\Gamma)$ is an $\SS(\fg^*)$-module. Moreover, since it sits inside 
the free $\SS(\fg^*)$-module \eqref{eq:1.4}, it is a torsion-free module. 
It is not, however, \emph{a priori} clear that it is itself a free module; 
this ``freeness'' issue is the first of the three questions we will take 
up below.

The second question concerns the number of generators of $H_G(\Gamma)$.
Let $\phi : V \to \RR$ be a Morse function, and, for 
every vertex, $p$, 
let $ind_{p}\phi$ be the number of vertices, $q$, adjacent to $p$, with 
$\phi(q) < \phi(p)$. This number, called \emph{the index of $p$} 
(and also denoted by $\sigma_p$), 
is the same as the number of oriented edges, $e$, with $i(e)=p$ and 
$\alpha_e(\xi) <0$. We define \emph{the $k$-th Betti number}, $b_k(\Gamma)$, 
to be the number of vertices with 
$ind_{p}\phi=k$~\footnote{In \cite{GZ1} and \cite{GZ2}, 
we used $b_{2k}(\Gamma)$ 
to denote this number; also, the $k$-th homogeneous component of 
$H_G(\Gamma)$ is denoted there by $H^{2k}(\Gamma, \alpha)$ and here 
by $H_G^k(\Gamma)$}.
One can show (\cite[Theorem 2.6]{GZ1}) that these numbers don't depend on 
$\xi$ or $\phi$, \emph{i.e.} are graph-theoretic invariants. Our question 
about the number of generators of $H_G(\Gamma)$ can be formulated as a 
conjecture:

\begin{conjecture}
The dimension of the $k$-th graded component of the ring
\begin{equation}
\label{eq:1.6}
H_G(\Gamma) \otimes_{\SS(\fg^*)} \CC
\end{equation}
is equal to $b_k(\Gamma)$.
\end{conjecture}

\begin{remark*}
It is clear that the orientation $o_{\xi}$ depends only on 
the connected component  of $\P$ in which $\xi$ sits. On the other hand it 
is clear that different components will give rise to different orientations. 
For instance, replacing $\xi$ with $-\xi$ reverses the orientation and 
changes a vertex of index $k$ into a vertex of index $d-k$. Therefore
\begin{equation}
\label{eq:poin}
b_k(\Gamma)= b_{d-k}(\Gamma).
\end{equation}
\end{remark*}

Finally, what do the generators of the ring $H_G(\Gamma)$ actually 
look like ?  For a vertex $p$ of $\Gamma$, let $F_p$ be the 
{\em flow-up of $p$}, {\em i.e.} the set of vertices of the 
oriented graph $(\Gamma, o_{\xi})$ that can be reached 
by a path along which the Morse function $\phi$ is strictly increasing.
($F_p$ is the graph theoretic analogue of the ``unstable manifold'' 
at a critical point, $p$, in ordinary Morse theory.)
We also define the {\em flow-down of} $p$ to be the set of vertices of the 
oriented graph $(\Gamma, o_{\xi})$ from where we can reach $p$ 
by a path along which the Morse function $\phi$ is strictly increasing.
Note that the flow-down of $p$ with 
respect to $o_{\xi}$ is the flow-up of $p$ with respect to
$o_{-\xi}$. Our conjecture is that,
just as in ordinary Morse theory, $H_G(\Gamma)$ is generated by 
cohomology classes dual to these ``unstable manifolds''. More explicitly:

\begin{conjecture}
There exists a set of independent generators
\begin{equation}
\label{eq:1.7}
\tau^{(p)} \in H_G^k(\Gamma), \qquad k = ind_{p}\phi
\end{equation}
with $supp \; \tau^{(p)} \subset F_p$. 
\end{conjecture}

Unfortunately, neither of these conjectures is true without some extra 
assumptions on $(\Gamma, \alpha)$. (An example of a pair $(\Gamma, \alpha)$ 
for which both conjectures fail to hold is described in 
\cite[Section~2.1]{GZ1}). 
Therefore, making a virtue of necessity, we will declare, by definition, 
that $(\Gamma, \alpha)$ satisfies \emph{the Morse package} if 
$H_G(\Gamma)$ is a free $\SS(\fg^*)$-module and if the two 
conjectures above are true.

We can now formulate the main result of this paper. Let 
$\fk$ be a linear subspace of $\fg$, $\fh=\fg/\fk$ and $\fh^*$ 
its dual, which we can identify with
\begin{equation}
\label{eq:none}
\fk^{\bot} = \{ \gamma \in \fg^* \; ; \; \gamma(\eta) = 0 
\mbox{ for all } \eta \in \fk \} \subset \fg^* .
\end{equation}
Define $\Gamma_{\fh^*}$ to be the subgraph of $\Gamma$ whose 
edges are all edges, $e$, of $\Gamma$ for which 
$\alpha_e \in \fh^*$.

\begin{theorem*}
\label{th:main1}
If $(\Gamma, \alpha)$ satisfies the Morse package, then, for every $\fh^*$ 
as above, $(\Gamma_{\fh^*}, \alpha)$ satisfies the Morse package.
Conversely, if $(\Gamma_{\fh^*}, \alpha)$ satisfies the Morse package 
for all $\fh^* \subset \fg^*$ of dimension 2, then $(\Gamma, \alpha)$ 
satisfies the Morse package.
\end{theorem*}

\begin{example}
{\bf The flag variety, $G/B$, $G=SL(n,\CC)$.} 
The graph $\Gamma$ is a 
Cayley graph associated with the permutation group, $S_n$. For this graph, 
$V = S_n$, and two permutations, 
$\sigma_k$, $k=1,2$, are adjacent 
in $\Gamma$ if and only if there exists a transposition 
$$\tau_{ij} \; : \;  i \leftrightarrow j \; , \qquad 1 \leq i < j \leq n$$
with $\sigma_2=\sigma_1 \tau_{ij}$. Moreover, if $e$ is the edge joining 
$\sigma$ to $\sigma\tau_{ij}$, the weight labeling $e$ is 
\begin{equation}
\label{eq:1.9}
\alpha_e = \begin{cases} 
\epsilon_j - \epsilon_i ,& \mbox{ if } \sigma(j)  > \sigma(i) \\
\epsilon_i - \epsilon_j ,& \mbox{ if } \sigma(j)  < \sigma(i),
\end{cases}
\end{equation}
$\epsilon_1,..,\epsilon_n$ 
being the standard basis vectors of the lattice $\ZZ^n$. The function 
\begin{equation}
\label{eq:1.10}
\phi : V \to \ZZ \; , \qquad \phi(\sigma) = \mbox{length}(\sigma)
\end{equation}
is a self-indexing Morse function on $V$; our theorem says that for 
$(\Gamma, \alpha)$ to satisfy the Morse package, it suffices to check that
the subgraphs associated with the subgroups $S_2$ and $S_3$ of $S_n$ 
satisfy the Morse package. This can be done more or less by inspection. 
For instance the graph associated with $S_3$ is the permutahedron 
(see Figure \ref{fig:perm}) 
and its ``Thom classes'' \eqref{eq:1.7} 
are given by simple monomials of 
degree 1, 2 and 3 in $\alpha_1, \alpha_2$ and $\alpha_1 + \alpha_2$,
with 
$\alpha_1= \epsilon_2 - \epsilon_1$ and $\alpha_2= \epsilon_3 - \epsilon_2$ 
(see the table below).
\begin{figure}[h]
\begin{center}
\includegraphics{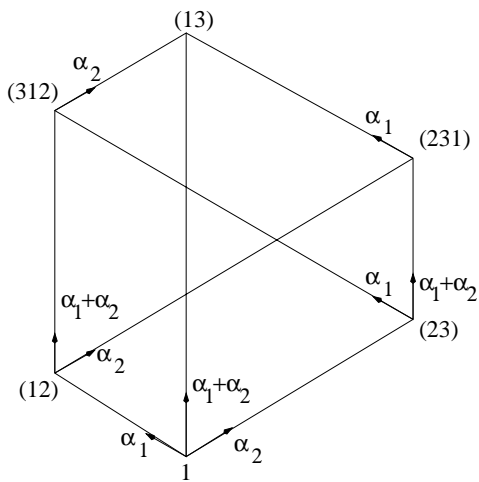}
\caption{Permutahedron}
\label{fig:perm}
\end{center}
\end{figure}
\begin{center}
\begin{tabular}{|c||c|c|c|c|c|c|} \hline
  & $\tau^1$ & $\tau^{(12)}$ & $\tau^{(23)}$ & $\tau^{(231)}$ & 
$\tau^{(312)}$ & $\tau^{(13)}$ \\ \hline \hline

 1 & 1 & 0 & 0 & 0 & 0 & 0 \\ \hline

 (12) & 1 & $-\alpha_1$ & 0 & 0 & 0 & 0 \\ \hline

 (23) & 1 & 0 & $-\alpha_2$ & 0 & 0 & 0 \\ \hline

 (231) & 1 & $-\alpha_1 - \alpha_2$ & $- \alpha_2$ &
 $\alpha_2 (\alpha_1 +\alpha_2)$ & 0 &
 0 \\ \hline

 (312) & 1 & $ -\alpha_1$ & $-\alpha_1 - \alpha_2$ &
 0 & $\alpha_1(\alpha_1 +\alpha_2)$ & 0 \\ \hline

 (13) & 1 & $-\alpha_1 -\alpha_2$ & 
 $- \alpha_1 - \alpha_2$ & $\alpha_2(\alpha_1 +\alpha_2)$ &
 $\alpha_1(\alpha_1 +\alpha_2)$ & $- \alpha_1 \alpha_2(\alpha_1 +\alpha_2)$
 \\ \hline
\end{tabular}
\end{center}
\end{example}

\smallskip

\begin{example}
{\bf The Grassmannian $Gr^k(\CC^n)$.}
Here the graph in question is the \emph{Johnson graph} 
(see \cite{BCN}). Its vertices are the $k$-element subsets, $S$, of the 
$n$-element set $\{ 1,..,n \}$, and two vertices, $S_k$, $k=1,2$, are 
adjacent if $S_1 \cap S_2$ is a $(k-1)$-element set, or, in other words, if
$S_1 - S_2$ and $S_2 - S_1$ are one-element sets. If $e$ is the oriented set 
joining $S_1$ to $S_2$ the weight labeling this edge is 
\begin{equation}
\label{eq:1.11}
\alpha_e = \alpha_i - \alpha_j,
\end{equation}
where $\{ i \} = S_1 - S_2$, $\{ j \} = S_2 - S_1$ and 
$\alpha_1,..,\alpha_n$ are the standard basis vectors of $\ZZ^n$. To define 
a Morse function on $V$ we first observe that there is a one-to-one
correspondence between the vertices of $\Gamma$ and Young diagrams 
(\cite{Fu}). Namely, 
let $S$ be a $k$-element subset of $\{ 1,..,n \}$, with elements 
$i_1 < i_2< .. < i_k$. Then we make correspond to $S$ the Young diagram, 
$\sigma_S$, with rows of length
\begin{equation}
\label{eq:1.12}
i_r - r \; , \qquad r = k,k-1,...,1.
\end{equation}
By means of this correspondence we can define a self-indexing Morse 
function, $\phi$, on $V$, by setting 
\begin{equation}
\label{eq:1.13}
\phi(S) = \mbox{ number of boxes in the diagram } \sigma_S.
\end{equation}

For the Johnson graph, our theorem says that to verify the Morse package 
for $\Gamma$, it suffices to verify the Morse package for the Johnson graph 
associated with the two-element subsets of $\{ 1,2,3 \}$, and just as in 
the previous example, this can be done more or less by inspection.
\end{example}

A few words about the contents of this paper: We mentioned above that the 
axioms A1-A3 can be translated into axioms on the axial function, $\alpha$.
In Section~\ref{sec:2} we will do this and also give more precise 
definitions of some of the concepts that we have mentioned in this 
introduction.
In Section~\ref{sec:3} we will prove the easy part  of our main theorem, 
the ``only if'' part. The hard part of the theorem, the ``if'' part, 
requires some algebraic results about Vandermonde matrices which we will 
discuss in the Appendix. Another ingredient in our proof is an 
\emph{integration} operation on graphs. Given an element, $f$, of the 
ring \eqref{eq:1.4}, we define its \emph{integral} 
over the graph to be the sum
\begin{equation}
\label{eq:1.14}
\int_{\Gamma} f = \sum_{p \in V} \frac{f_p}{\prod \alpha_e},
\end{equation}
the denominator in the $p^{th}$ 
summand being the product over all oriented 
edges, $e$, with  $i(e)=p$. The expression on the right appears to belong 
to the quotient field of $\SS(\fg^*)$; however, if 
$f \in H_G(\Gamma)$, one can show that this sum is in $\SS(\fg^*)$. 
Moreover, $H_G(\Gamma)$ has the following duality property:
\begin{proposition}
\label{prop:dual}
If $f \in \mbox{Maps}(V, \SS(\fg^*))$ and 
$$\int_{\Gamma} fh  \; \in \; \SS(\fg^*)$$
for all $h \in H_G(\Gamma)$, then $f \in H_G(\Gamma)$.
\end{proposition}

These results are old results of ours from \cite{GZ1}, but to prove the 
main theorem, we will need to generalize them to \emph{hypergraphs}. 
(A hypergraph is, like a graph, an object consisting of vertices, 
hyperedges and an incidence relation which tells one when a vertex is incident 
to a hyperedge. For hypergraphs, however, more than two vertices can be 
incident to a hyperedge.)

A third ingredient in our proof is a ``cross-section'' construction for
graphs which we already made extensive use of in \cite{GZ2}. A real number, 
$c$, is a \emph{critical value} of the Morse function $\phi$ if 
$c=\phi(p)$ for some $p \in V$; otherwise, it is a 
\emph{regular value}. Let $c$ be a regular value and let $V_c$ be the set 
of oriented edges, $e$, of $\Gamma$, which intersect the level set 
$\phi = c$, in the sense that 
$\phi(p) < c < \phi(q),$
$p=i(e)$ and $q=t(e)$ being the vertices of $e$. 
Then $V_c$ is the vertex set of a 
hypergraph, $\Gamma_c$, whose edges are the intersections of the level set 
$\phi=c$ with the connected components of the graphs $\Gamma_{\fh^*}$, 
$\fh^*$ being a dimension 2 subspace of $\fg^*$. 

The trickiest and most subtle part of our proof is figuring out how to define 
the cohomology ring of this hypergraph. One can equip $\Gamma_c$ with an 
analogue of the axial function, $\alpha$, (see Section \ref{sec:localiz}) and 
mimic the definition \eqref{eq:1.5} of $H_G(\Gamma)$; however, this 
turns out \emph{not} to be the right approach. Trial and error show that 
the right approach is to mimic the definition of $H_G(\Gamma)$ 
``by duality'' in Proposition \ref{prop:dual}. (See Section \ref{sec:5} 
for details). This approach, however, has one apparent short-coming: 
The equivariant cohomology theory of graphs that we described above 
is a \emph{local} theory in the sense that a prospective cohomology 
class, $f \in \mbox{Maps }(V, \SS(\fg^*))$, belongs to $H_G(\Gamma)$ 
if and only if it satisfies condition \eqref{eq:1.5} at each edge. If we 
define hypergraph cohomology by duality it is not at all clear that it 
will be a local theory in this sense. The localization theorem which 
we prove in Section \ref{sec:localiz} shows, however, that it is.

Denoting the cohomology ring of $\Gamma_c$ 
by $H(\Gamma_c, \alpha_c)$, we show in Section~\ref{sec:6}
how the proof of the main theorem can be reduced to the following problem: 
How does the structure of $H(\Gamma_c, \alpha_c)$ change as one passes 
through a critical value of $\phi$ ? 
If $c$ and $c'$ are regular values of $\phi$ 
and there is just one critical point, $p \in V_{\Gamma}$, with 
$c < \phi(p) < c'$, how are $H(\Gamma_c, \alpha_c)$ and 
$H(\Gamma_{c'}, \alpha_{c'})$ related ?
The main technical result of this paper answers this question by showing 
that the elements of $H(\Gamma_{c'}, \alpha_{c'})$ are expressible in terms 
of the elements of $H(\Gamma_c, \alpha_c)$ via an ``integral'' equation of 
Vandermonde type.

\smallskip

\begin{center}
\section{The Main Result}
\label{sec:2}
\end{center}

\smallskip

Let $\Gamma$ be a regular $d$-valent graph, $V$ its set of vertices, 
and, for each vertex $p \in V$, 
let  $E_p$ be the set of oriented  edges with initial vertex $p$. 
Define the ``tangent bundle'' of $\Gamma$ to be ($E_{\Gamma}$, $\pi$), 
with projection $\pi : E_{\Gamma} \to V_{\Gamma}$, $\pi(e)=i(e)$.
Let $\fg^*$ be an $n$-dimensional real vector space.

\begin{definition}
\label{def:abs1skel}
An \emph{abstract one-skeleton} 
is a pair $(\Gamma, \alpha)$ consisting of a $d$-valent graph, $\Gamma$, 
and a function 
$$\alpha : E_{\Gamma} \to \fg^*$$
(called {\it axial function}), satisfying the axioms:

\begin{enumerate}[{Axiom A}1]

\item :
For every $p \in V$, the vectors 
$\{ \alpha_e \; ; \; e \in E_p \}$ 
are pairwise linearly independent.

\item : 
If $e \in E_{\Gamma}$ then
\begin{equation*}
\alpha_{\bar{e}} = - \alpha_{e}.
\end{equation*}

\item :
Let $e \in E_{\Gamma}$ with $p=i(e)$ and $p'=t(e)$. 
Let $E_p = \{e_i; i=1,..,d \}$ and $E_{p'} = \{ e_i'; i=1,..,d\}$. 
Then  we can order $E_p$ and $E_{p'}$ such that $e_d=e$, $e_d'=\bar{e}$ 
and
\begin{equation*}
\alpha_{e_i'} = \alpha_{e_i} + c_{i,e} \alpha_e \; , 
\qquad \mbox{with } c_{i,e} \in \RR.
\end{equation*}
\end{enumerate}
\end{definition}

If $\fh^*$ is a linear subspace of $\fg^*$, we define $\Gamma_{\fh^*}$ 
to be the subgraph of $\Gamma$ whose edges are all edges $e$ for which 
$\alpha_e \in \fh^*$. Then $(\Gamma_{\fh^*}, \alpha)$ is an abstract 
one-skeleton, whose axial function takes values in $\fh^*$. 
Note that if $\fk$ is a subspace of $\fg$ and 
$\fh = \fg/\fk$ then 
$\fh^*$ is canonically a subspace of $\fg^*$ and  
$$\Gamma_{\fh^*} \simeq \Gamma_{(\fg/\fk)^*}.$$

\begin{definition}
A vector $\xi \in \fg$ is called {\em generic} if
for every vertex $p \in V$ and every four edges 
$e_1,e_2,e_3$ and $e_4$ in  $E_p$ such that $e_1 \neq e_2$, $e_3 \neq e_4$
and $(e_1,e_2) \neq (e_3,e_4)$, we have
\begin{equation*}
\frac{\alpha_{e_1}}{\alpha_{e_1}(\xi)} - 
\frac{\alpha_{e_2}}{\alpha_{e_2}(\xi)}
 \neq 
\frac{\alpha_{e_3}}{\alpha_{e_3}(\xi)} - 
\frac{\alpha_{e_4}}{\alpha_{e_4}(\xi)}.
\end{equation*}
\end{definition}
The set of generic elements, $\P_0$, is dense in $\fg$. 

From now on, we assume that $(\Gamma, \alpha)$ is a $\xi$-acyclic
one-skeleton for some generic and polarizing $\xi$ and that 
$\phi : V \to \RR$ is $\xi$-compatible. 
Without loss of generality, we can assume that $\phi$ is injective.
Let $H(\Gamma, \alpha)$ be the $\SS(\fg^*)$-module defined by 
\eqref{eq:1.5}.

\begin{definition}
An element $\tau^{(p)}$ of $H(\Gamma,\alpha)$ is called 
\emph{a generating class at} $p$ if it satisfies the conditions:
\begin{enumerate}
\item $\tau^{(p)}$ is supported on the flow-up at $p$, $F_p$
\item $\tau^{(p)}(p)= \prod' \alpha_{e}$
\end{enumerate}
the product in the second condition being over the edges 
$e \in E_p$ with  $\alpha_{e}(\xi)<0$. A set $\{ \tau^{(p)} \}_{p\in V}$, 
containing one generating class for each vertex,
is called \emph{a generating family}.
\end{definition}

\begin{theorem}
\label{th:equivalence}
The one skeleton $(\Gamma,\alpha)$ admits a generating family if and only if
for every $m \geq 0$,
\begin{equation}
\label{eq:dimens}
\dim H^m(\Gamma,\alpha) = \sum_{k=0}^{d} b_{k}(\Gamma) \lambda_{m-k,n},
\end{equation}
where
\begin{equation*}
\lambda_j = \lambda_{j,n} = 
\begin{cases}
\dim \SS^j(\fg^*) & , \mbox{ if } j \geq 0 \\
0               & , \mbox{ if } j < 0. 
\end{cases}
\end{equation*}
\end{theorem}

\begin{proof}
The equivalence follows from Theorems 2.4.2 and 2.4.4 of \cite{GZ2}.
\end{proof}

\begin{remark*}
A generating family is a basis of $H(\Gamma,\alpha)$ as an 
$\SS(\fg^*)$-module. So, if a generating family exists, $H(\Gamma,\alpha)$ 
is a free $\SS(\fg^*)$-module.
\end{remark*}

The main result of this paper is the following theorem:

\begin{theorem}
\label{th:mainres}
The one-skeleton $(\Gamma, \alpha)$ admits a generating family
if and only if for every two-dimensional subspace $\fh^* \subset \fg^*$, 
every connected component of the induced one-skeleton 
$(\Gamma_{\fh^*}, \alpha)$ admits a generating family.
\end{theorem}

One also has the following sharpening of the ``if'' part of this theorem:

\begin{theorem}
\label{th:sharpening}
Suppose that for every $q \in F_p - \{p \}$, the index of $q$ is strictly 
greater that the index of $p$. Then the class $\tau^{(p)}$ is unique.
\end{theorem}

In view of Theorem \ref{th:equivalence}, to prove 
Theorem \ref{th:mainres}, it suffices to show that 
for every $m$, the dimension of $H^m(\Gamma,\alpha)$ is given by
\eqref{eq:dimens}. 
To compute this dimension, we first realize $(\Gamma, \alpha)$ as a 
cross-section of a product one-skeleton $\Gamma^{\flat} = \Gamma \times L$ 
using a process akin to ``symplectic cutting'' (\cite{Le}).
To each cross-section we attach a module and we show that the module 
we associate to a cross-section ``isomorphic'' to $(\Gamma, \alpha)$ 
is the same as the cohomology ring of $(\Gamma, \alpha)$. 
When the level is just above the lowest vertex, the cross-section
is very simple and we determine the associated module directly.
Then we determine how this module
changes when we pass from one cross-section to another, going over a vertex.
Finally, we compute the dimension of $H^m(\Gamma,\alpha)$ 
by starting with the dimension of the first cross-section and 
counting the changes that occur until we have reached the cross-section of 
$\Gamma^{\flat}$ that is the same as $(\Gamma, \alpha)$.

\smallskip

\begin{center}
\section{The ``only if'' part of Theorem \ref{th:mainres}}
\label{sec:3}
\end{center}

\smallskip

Let $\fh^*$ be a two-dimensional subspace of $\fg^*$,
$\fg$ be the dual of $\fg^*$ and $\fk$ be the annihilator of 
$\fh^*$ in $\fg$. Choose $\fh$ to be a complementary subspace to 
$\fk$ in $\fg$; then we can identify $\fh^*$ with the vector space dual to 
$\fh$. 

Now let $p$ be a vertex of $(\Gamma_{\fh^*}, \alpha)$ and let $\tau$ be 
the Thom class, in $H(\Gamma,\alpha)$, associated with the 
flow-up of $p$. Thus
\begin{equation}
\label{eq:taupe}
\tau_p = \prod_{i=1}^r \alpha_{e_i},
\end{equation}
where $r$ is the index of $p$ and $e_1$,...,$e_r$ are the 
``downward-pointing'' edges at $p$, \emph{i.e.} $\alpha_{e_i}(\xi) <0$. 
Suppose that $\alpha_{e_i} \in \fh^*$ for
$i=1,..,s$ and $\alpha_{e_i} \not \in \fh^*$ for $i=s+1,..,r$. Choose
$\xi \in \fk$ such that $\alpha_{e_i}(\xi) \neq 0$ for 
$i=s+1,..,r$. Given $f \in \SS(\fg^*)$, let $f^{\#}$ be the restriction of
this polynomial function on $\fg$ to the two-dimensional affine subspace 
$$\{ \xi \} \times \fh \subset \fk \times \fh = \fk \oplus \fh = \fg.$$
Then we can think of $f^{\#}$ as an element of $\SS(\fh^*)$. 

\begin{lemma}
\label{lem:vi1}
Let $V_{\fh^*}$ be the set of vertices of $\Gamma_{\fh^*}$. For 
$h \in H(\Gamma,\alpha)$, 
let $h^{\#} : V_{\fh^*} \to \SS(\fh^*)$ be defined by
$$(h^{\#})(p) = (h(p))^{\#}.$$
Then $h^{\#} \in H(\Gamma_{\fh^*},\alpha)$.
\end{lemma}

\begin{proof}
If $p$ and $q$ are in $V_{\fh^*}$ 
and are joined by an edge $e$ of $\Gamma_{\fh^*}$, then 
$$h(p) - h(q) = f \alpha_e$$
for some $f \in \SS(\fg^*)$, so that
$$h^{\#}(p) - h^{\#}(q) = f^{\#} (\alpha_e)^{\#}.$$
However, since $\alpha_e \in \fh^*$, we have $(\alpha_e)^{\#}=\alpha_e$.
\end{proof}

Consider, in particular, $\tau^{\#}$. By \eqref{eq:taupe},
$$\tau_p^{\#} = \Bigl( \prod_{i=s+1}^r \alpha_{e_i}^{\#} \Bigr)
\prod_{i=1}^s \alpha_{e_i}.$$
Moreover, for $s+1 \leq i \leq r$,
$$\alpha_{e_i}^{\#} = c_i + \beta_{e_i},$$
where $c_i = \alpha_{e_i}(\xi) \neq 0$ and $\beta_{e_i} \in \fh^*$.
Thus, if $c$ is the product of the $c_i$'s, 
$$\tau_p^{\#} =  c \prod_{i=1}^s \alpha_{e_i} + g,  \;\; \mbox{ with } 
\;g \in \bigoplus_{j \geq s+1}\SS^j(\fh^*).$$
Let $\tau_{\fh^*}^{(p)}$ be the homogeneous component of degree $s$ of the 
cohomology class $c^{-1}\tau^{\#}$. Then 
$$ \tau_{\fh^*}^{(p)}(p) = \prod_{i=1}^s \alpha_{e_i}$$
and, since $\tau^{\#}$ is supported on the flow-up of $p$ in $\Gamma_{\fh^*}$,
so is $\tau_{\fh^*}^{(p)}$.

\smallskip

\begin{center}
\section{Cross-sections}
\label{sec:5}
\end{center}

\smallskip

We will say that an axial function, $\alpha$, is 3-independent if, 
for every $p \in V$ and every triple of edges, $e_i \in E_p$, $i=1,2,3$, 
the vectors $\alpha_{e_i}$ are linearly independent. In \cite{GZ2} the
``if'' part of Theorem \ref{th:mainres} was proved modulo this 3-independence 
assumption. We will \emph{not} make this assumption here, however, we will 
make use of a key ingredient in our earlier proof: Fixing a regular value 
of $\phi$, $c \in \RR - \phi(V)$, we will define the $c$-cross-section, 
$V_c$, of $\Gamma$, to be the set of oriented edges, $e$, of 
$\Gamma$, with end points $p=i(e)$ and $q=t(e)$, for which
\begin{equation}
\label{eq:new1}
\phi(p) < c < \phi(q).
\end{equation}
In \cite{GZ2} we showed that $V_c$ is the set of vertices of a graph, 
$\Gamma_c$, whose edges are defined as follows: For every two dimensional 
subspace, $\fh^*$, of $\fg^*$, and every connected component, 
$\Gamma_{\fh^*}^0$, of $\Gamma_{\fh^*}$, let $V_c(\Gamma_{\fh^*}^0)$ be the 
edges of $\Gamma_{\fh^*}$ satisfying \eqref{eq:new1}. If $\alpha$ is 
3-independent and if the ``only if'' part of Theorem \ref{th:mainres} holds, 
these subsets of $V_c$ are of cardinality 2 and are by 
definition the edges of $\Gamma_c$. 

If we drop the 3-independence assumption we can still define these sets to 
be the \emph{hyperedges} of $\Gamma_c$, but the object we get will no longer 
be a graph but, instead, a \emph{hypergraph}. Nonetheless many of the results 
of \cite{GZ2} will continue to hold. We will describe some of these results 
in this section.

\smallskip

Let $\fg^*_{\xi}$ be the annihilator of $\xi$ in $\fg^*$, $\SS(\fg_{\xi}^*)$ 
be the symmetric algebra of $\fg^*_{\xi}$ and $Q(\fg^*_{\xi})$ be the 
quotient field of $\SS(\fg_{\xi}^*)$. We define a morphism of graded rings
\begin{equation}
\label{eq:kir} 
 \K_c : H(\Gamma,\alpha) \to \mbox{Maps}(V_c, \SS(\fg_{\xi}^*))
\end{equation}
as follows: For $e \in V_c$ let $p$ and $q$ be the vertices
of $e$. Since $\alpha_{e}(\xi) \ne 0$,
the projection $\fg^* \to \fg^*_e$ maps $\fg^*_{\xi}$
bijectively onto $\fg^*_e$, so one
has a composite map
\begin{displaymath}
\fg^* \to \fg^*_e \leftrightarrow \fg^*_{\xi}
\end{displaymath}
and hence an induced morphism of graded rings:
\begin{displaymath}
\SS(\fg^*) \to \SS(\fg^*_e) \leftrightarrow \SS(\fg^*_{\xi}) \, .
\end{displaymath}
If $f$ is in $H(\Gamma,\alpha)$, the images of $f_p$ and 
$f_q$ in $\SS(\fg^*_e)$
are the same by \eqref{eq:rhoe} and hence so are their images in
$\SS(\fg^*_{\xi})$.  We define $\K_c(f)(e)$ to be this common image
and call the map $\K_c$ {\em the Kirwan map}.

Let $\{ x,y_1,...,y_{n-1} \}$ be a basis of $\fg^*$ such that 
$x(\xi)=1$ and $\{ y_1,...,y_{n-1} \}$ is a basis of $\fg_{\xi}^*$.
Let $\alpha = \alpha(\xi)(x-\beta(y)) \in \fg^*$ such that
$\alpha(\xi) \ne 0$. Then the map 
$$\rho_{\alpha}: \SS(\fg^*) \rightarrow \SS(\fg_{\xi}^*),$$
given via the 
identification $\fg_{\xi}^* \simeq \fg^*/ \RR\alpha$, will send 
$x$ to $x-\alpha/\alpha(\xi)\in \fg_{\xi}^*$ and $y_j$ to $y_j$.
Therefore $\rho_{\alpha}$ will send a polynomial 
$P(x,y) \in \SS(\fg^*)$ to the polynomial 
$P(x-\alpha/\alpha(\xi),y) = P(\beta(y),y)\in \SS(\fg_{\xi}^*)$.

Hence, if $\alpha_{e}= m_e(x-\beta_e(y))$, then 
\begin{equation*}
\K_c(f)(e) = f_p(\beta_e(y),y) \in \SS(\fg_{\xi}^*).
\end{equation*}

Let $e \in V_c$, with endpoints $p=i(e)$ and $q=t(e)$, such that
$\phi(p) < c < \phi(q)$. 
If $e_i$, $i=1,..,d-1$ are the other edges issuing from $p$ and 
$e_i'$, $i=1,..,d-1$ are the other edges issuing from $q$, 
we define the Thom class of $e$ in $\Gamma$, 
$\tau_e \in H(\Gamma,\alpha)$, by
\begin{equation*}
\tau_e(v) = 
\begin{cases} 
\prod \alpha_{e_i} & \text{ if  $v=p$} \\
\prod \alpha_{e_i'} & \text{ if $ v=q$} \\
0 &  \text{ otherwise}.
\end{cases}
\end{equation*}
Let
\begin{equation*}
\alpha_i^{\#} = \K_c(\alpha_i) = 
\alpha_i -\frac{\alpha_i(\xi)}{\alpha_e(\xi)}\alpha_e 
\end{equation*}
and
\begin{equation}
\label{eq:deltae}
\delta_e = (m_e \K_c(\tau_e))^{-1} =
(m_e \prod_{i=1}^{d-1} \alpha_i^{\#})^{-1} 
\in Q(\fg_{\xi}^*)
\end{equation}
where $m_e = \alpha_{e}(\xi) > 0$.

We now define the integral operation
\begin{equation*}
\int_{\Gamma_c} : \mbox{Maps}(V_c, Q(\fg_{\xi}^*)) \to Q(\fg_{\xi}^*)
\end{equation*}
as integration with respect to the discrete measure on $\Gamma_c$
with density $\delta_e$ at $e$, {\it i.e.}
\begin{equation}
\label{eq:intgc}
\int_{\Gamma_c} f = \sum_{e\in V_c} \delta_e f(e) = 
\sum_{e\in V_c} \frac{f(e)}
{m_e \prod_{i=1}^{d-1} \alpha_i^{\#}}.
\end{equation}

One of the main results in \cite{GZ1} (Theorem 2.5) is :

\begin{theorem}
\label{th:kir}
If $f \in H(\Gamma,\alpha)$ then
\begin{equation*} 
\int_{\Gamma_c} \K_c(f) \in \SS(\fg_{\xi}^*).
\end{equation*}
\end{theorem}

We will use this property to associate an $\SS(\fg_{\xi}^*)$-module 
to the cross-section.

\begin{definition}
\label{def:cohcross}
We denote by $H(\Gamma_c, \alpha_c)$ the set of all maps
$f : V_c \to Q(\fg_{\xi}^*)$ with the property that
\begin{equation*}
\int_{\Gamma_c} f \K_c(h) \in \SS(\fg_{\xi}^*)
\end{equation*}
for all $h \in H(\Gamma,\alpha)$.
\end{definition}

\noindent{\bf Remarks:}

\begin{enumerate}
\item Let $f \in H(\Gamma_c, \alpha_c)$. Then for 
every edge $e \in \Gamma_c$, 
$$f(e)= m_e \int_{\Gamma_c}  f \K_c(\tau_e) \in \SS(\fg_{\xi}^*).$$
Therefore 
$H(\Gamma_c, \alpha_c) \subset \mbox{Maps}(V_c, \SS(\fg_{\xi}^*))$.

\item
Since the Kirwan map \eqref{eq:kir} is a morphism of rings,
from Theorem \ref{th:kir} follows that the image of the Kirwan map is 
included in $H(\Gamma_c, \alpha_c)$.
\end{enumerate}

\begin{example}
Let $p$ be the vertex of $\Gamma$ on which $\phi$ attains its minimum 
and let $c \in \RR$ be a regular value such that $p$ is the only vertex 
below the $c$-level. Then the $c$-cross section consists of the oriented 
edges with initial vertex at $p$.

For such an edge, $e_i$, let
\begin{equation*}
\alpha_{e_i} = m_i(x -\beta_i(y)),
\end{equation*}
with $m_i = \alpha_{e_i}(\xi)$ and $\beta_i(y) \in \SS^1(\fg_{\xi}^*)$.
Consider $\tau : V_c \to \fg_{\xi}^*$, $\tau(e_i) = \beta_i(y)$.

Let $g : V_c \to \SS(\fg_{\xi}^*)$.
For $h \in H(\Gamma,\alpha)$, let 
$P=h(p)\in \SS(\fg^*) \simeq \SS(\fg_{\xi}^*)[x]$.
Then
\begin{equation*}
\int_{\Gamma_c} g\K_c(h) = (\prod_{i=1}^d m_i)^{-1}
\int_{V_c} g P(\tau) \quad ,
\end{equation*}
where the integral on the right hand side is defined as in 
\eqref{eq:intonfinite} in the Appendix.
Therefore, by Lemma \ref{lem:14.1} in the Appendix, 
$g$ is an element of $H(\Gamma_c, \alpha_c)$ if and 
only if there exist $g_0,...,g_{d-1} \in \SS(\fg_{\xi}^*)$ such that
\begin{equation}
\label{eq:genelem}
g=\sum_{k=0}^{d-1} g_k\tau^k.
\end{equation}
We conclude that 
\begin{equation*}
\mbox{dim}H^{m}(\Gamma_c,\alpha_c) = \sum_{k=0}^{d-1} \lambda_{m-k,n}.
\end{equation*}

Moreover, let $g \in  H(\Gamma_c, \alpha_c)$ be given by 
\eqref{eq:genelem} and consider 
\begin{equation*}
h = \sum_{k=0}^{d-1} g_k x^k \in \SS(\fg_{\xi}^*)[x] \simeq 
\SS(\fg^*) \subset H(\Gamma,\alpha).
\end{equation*}
Then $g = \K_c(h)$; this proves that for this cross-section, 
the Kirwan map $\K_c$ is surjective.
\end{example}

\smallskip

\begin{center}
\section{The localization theorem for $H(\Gamma_c,\alpha_c)$}
\label{sec:localiz}
\end{center}

At the beginning of Section \ref{sec:5} we pointed out that the set $V_c$ 
can be made into the set of vertices of a hypergraph, $\Gamma_c$. In this 
section we will describe the structure of this hypergraph in more detail 
and give an alternative description of $H(\Gamma_c, \alpha_c)$, which is 
more in the spirit of the description \eqref{eq:1.5} of $H(\Gamma,\alpha)$.

Recall that an edge of $\Gamma_c$ is a subset of $V_c$ of the form
\begin{equation}
\label{eq:new5.1}
E= V_c(\Gamma_{\fh^*}^0),
\end{equation}
where $\fh^*$ is a 2-dimensional subspace of $\fg^*$, spanned by axial 
vectors $\alpha_{e_1}$ and $\alpha_{e_2}$, $e_i \in E_p$, and 
$\Gamma_{\fh^*}^0$ is the connected component of $\Gamma_{\fh^*}$ 
containing $p$. If one assigns to the edge \eqref{eq:new5.1} the 
multiplicity
$$\mu_E = (\mbox{ valence of } \Gamma_{\fh^*}^0 ) -1 \; ,  $$
then $\Gamma_c$ becomes a $(d-1)$-valent hypergraph: each vertex is the 
point of intersection of $d-1$ edges, counting multiplicities. Let us label 
each of the edges, $E$, by a non-zero vector, $\alpha_E$, in the 
one-dimensional vector space $\fh^* \cap \fg_{\xi}^*$. This labeling is not 
unique since $\alpha_E$ is only defined up to a scalar; however, one can 
show that it satisfies axioms similar to the axioms A1-A3 of definition 
\ref{def:abs1skel} (see \cite[Section 2.1]{GZ2}). Now fix a vector 
$\gamma \in \fg_{\xi}^* - 0 $ and let $\Gamma_c^{\gamma}$ be the subgraph 
of $\Gamma_c$ consisting of the edges, $E$, of $\Gamma_c$, with 
$\alpha_E \in \RR \gamma$.

\begin{lemma}
\label{lem:new5.1}
The hypergraph $\Gamma_c^{\gamma}$ is totally disconnected: no two distinct 
edges, $E_1$ and $E_2$, of $\Gamma_c^{\gamma}$, contain a common vertex.
\end{lemma}

\begin{proof} Let $v$ be a vertex of $\Gamma_{c}^{\gamma}$, corresponding to 
an edge $e$ of $\Gamma$. Let $E$ be an edge of $\Gamma_{c}^{\gamma}$ incident 
to $v$ and let $\fh^*$ be the two-dimensional subspace of $\fg^*$ such that 
$ E= V_c(\Gamma_{\fh^*}^0)$. Then $\fh^*$ is generated by $\gamma$ and 
$\alpha_e$, therefore is uniquely determined.
\end{proof}

As in \eqref{eq:new5.1}, let $\fh^*$ be the 2-dimensional subspace of $\fg^*$
spanned by vectors $\alpha_{e_1}, \alpha_{e_2}$, $e_1, e_2 \in E_p$ and let 
$\Gamma_{\fh^*}^0$ be the connected component of $\Gamma_{\fh^*}$ 
containing $p$. 

\begin{definition}
\label{def:new1}
The \emph{Thom class} of $\Gamma_{\fh^*}^0$,
$$\tau = \tau(\Gamma_{\fh^*}^0) \in H(\Gamma,\alpha) \; , $$
is the class which assigns to each vertex $q$ of $\Gamma_{\fh^*}^0$ the 
product
$$\prod_{e_i \in E_q -E_q(\Gamma_{\fh^*}^0)} \alpha_{e_i}$$
and assigns zero to every other vertex of $\Gamma$.
\end{definition}

The existence of such a class enables one to define an injection
$$H(\Gamma_{\fh^*},\alpha) \to H(\Gamma,\alpha) \; , \qquad f \mapsto f\tau.$$
Moreover, from $\tau$ one gets an important identity relating integration 
over $\Gamma_c$ and integration over the edges \eqref{eq:new5.1} of 
$\Gamma_c$:

\begin{lemma}
\label{lem:new5.2}
Given $g \in H(\Gamma_{\fh^*}^0,\alpha)$ and $f : V_c \to \SS(\fg_{\xi}^*)$,
one has
\begin{equation}
\label{eq:new5.2}
\int_E f_E \K_c(g) = \int_{\Gamma_c} f \K_c(\tau g),
\end{equation}
where $f_E$ is the restriction of $f$ to $E$ and the Kirwan map 
$\K_c$ on the left hand side is the map
$$\K_c : H(\Gamma_{\fh^*}^0,\alpha) \to 
H((\Gamma_{\fh^*}^0)_c, \alpha_c).$$
\end{lemma}

We will use Lemma \ref{lem:new5.2} to prove:

\begin{theorem}
\label{th:new5}
A map $f : V_c \to \SS(\fg_{\xi}^*)$ is in $H(\Gamma_c, \alpha_c)$ 
if and only if its restriction to every edge $E$ of the hypergraph 
$\Gamma_c$ is in $H(E,\alpha_c)$.
\end{theorem}

\begin{proof}
Assume first that $f \in H(\Gamma_c, \alpha_c)$ and let $E$ be an edge of
$\Gamma_c$, given by \eqref{eq:new5.1}. Denote by $f_E$ the restriction 
of $f$ to $E$. For every $g \in  H(\Gamma_{\fh^*}^0,\alpha)$ we have 
$$\int_E f_E \K_c(g) = \int_{\Gamma_c} f \K_c(\tau g) \in \SS(\fg_{\xi}^*),$$
hence $f_E \in H(E,\alpha_c)$.

Conversely, assume that $f_E \in H(E,\alpha_c)$ for every edge $E$ of 
$\Gamma_c$. Let $g \in H(\Gamma,\alpha)$. We need to show that 
\begin{equation}
\label{eq:zero}
\int_{\Gamma_c} f \K_c(g) \in \SS(\fg_{\xi}^*).
\end{equation}
where, according to \eqref{eq:intgc},
\begin{equation}
\label{eq:prima}
\int_{\Gamma_c} f \K_c(g) = \sum_{e \in V_c} \delta(e) f(e) \K_c(g)(e).
\end{equation}
Let $\M \in \fg_{\xi}^*$ be the set of all linear factors that appear in 
the denominators of $\delta(e)$ (see \eqref{eq:deltae}), for all edges $e$
that intersect the $c$-level. For $\gamma \in \M$, let 
$\M_{\gamma}$ be the multiplicative system in $\SS(\fg_{\xi}^*)$ 
generated by elements of $\M$ that are not multiples of $\gamma$ and let 
$\SS(\fg_{\xi}^*)_{\gamma}$ be the corresponding localized ring.
To show \eqref{eq:zero} it suffices to show that
\begin{equation}
\label{eq:doi}
\int_{\Gamma_c} f \K_c(g) \in \SS(\fg_{\xi}^*)_{\gamma}
\end{equation}
for every $\gamma \in \M$. 

Fix $\gamma \in \M$ and let $\Gamma_c^{\gamma}$ be the corresponding 
subgraph of $\Gamma_c$. Let $e$ be a vertex of $\Gamma_c^{\gamma}$ 
corresponding to an edge of $\Gamma$, and let $E=V_c(\Gamma_{\fh^*}^0)$ 
be an edge of $\Gamma_c^{\gamma}$ incident to $e$, where $\fh^*$ is a 
two-dimensional subspace of $\fg^*$. Then $e$ is an edge in both $\Gamma$
and $\Gamma_{\fh^*}^0$; the density associated to $e$ when we integrate 
over $\Gamma_c$ is given by \eqref{eq:deltae} and the density associated 
to $e$ when we integrate over $E$ is 
\begin{equation}
\label{eq:deltae'} 
\delta_e' = (m_e\K_c(\tau_{e}'))^{-1},
\end{equation}
where $\tau_{e}'$ is the Thom class of $e$ in $\Gamma_{\fh^*}^0$. 
If $\tau \in H(\Gamma,\alpha)$ is the Thom class of $\Gamma_{\fh^*}^0$ and
\begin{equation}
\label{eq:ne}
\fn_E(e) =  \K_c(\tau)(e) \in \MM_{\gamma} \; , 
\qquad \mbox{ for all } e \in E
\end{equation}
then the two densities are related by 
\begin{equation}
\label{eq:2dens}
\delta_e = \frac{\delta_e'}{\fn_E(e)}.
\end{equation}
Then \eqref{eq:prima} differs from
\begin{equation}
\label{eq:trei}
\sum_{E} \sum_{e \in E} 
\delta_e f(e) \K_c(g)(e) =
\sum_{E} \sum_{e \in E} 
\frac{\delta_e' f(e) \K_c(g)(e)}{\fn_E(e)}
\end{equation}
by a term in $\SS(\fg_{\xi}^*)_{\gamma}$ (the first sums on each side are over
all edges of $\Gamma_c^{\gamma}$).

But 
\begin{equation}
\label{eq:patru}
\sum_{e \in E} \frac{\delta_e' f(e) \K_c(g)(e)}{\fn_E(e)} =
\frac{1}{\prod_{e \in E} \fn_E(e)} 
\sum_{e \in E} \Bigl( \delta_e' f(e) \K_c(g)(e) 
\prod_{e' \neq e} \fn_E(e') \Bigr),
\end{equation}
and, according to Lemma \ref{lem:2} of the Appendix, there exists a 
polynomial $P_E \in \SS(\fg_{\xi}^*)[X]$ such that 
\begin{equation}
\prod_{e' \neq e} \fn_E(e') = P_E(\fn_E(e));
\end{equation}
then, since $\K_c$ is a morphism of rings, \eqref{eq:trei} becomes 
\begin{equation}
\sum_{E}\frac{1}{\prod_{e \in E} \fn_E(e)} 
\int_E f_E \K_c(gP_E(\tau)) \in \SS(\fg_{\xi}^*)_{\gamma},
\end{equation}
which concludes the proof.
\end{proof}

\smallskip

\begin{center}
\section{Cutting}
\label{sec:cutting}
\end{center}

\smallskip

Let $L$ be a one-valent graph, with two vertices, labeled 0 and 1, and one
edge connecting them. Consider the product graph 
$\Gamma^{\flat} = \Gamma \times L$, with the set of vertices $V^{\flat}$; 
we define an axial function
$$\alpha^{\flat} : E_{\Gamma^{\flat}} \to (\fg \oplus \RR)^* 
\simeq \fg^* \oplus \RR^* \; \; ,$$  
by
\begin{align*}
\alpha^{\flat}_{(p,0)(q,0)} & = \alpha_{pq}, 
& \mbox{ if } pq \in E_{\Gamma} \\
\alpha^{\flat}_{(p,1)(q,1))} & = \alpha_{pq}, 
&\mbox{ if } pq \in E_{\Gamma} \\
\alpha^{\flat}_{(p,0)(p,1)} & =  {\bf 1}\in \RR^* \; ; &
\end{align*}
(this would correspond to the product action of $G \times S^1$ on 
$\Gamma^{\flat}=\Gamma \times L$).

Then the pair $(\Gamma^{\flat}, \alpha^{\flat})$ is a one-skeleton, 
acyclic with respect to the generic vector 
$(\xi,1) \in \fg \oplus \RR$. Choose 
$a > \phi_{max} - \phi_{min} > 0$ and define
$\Phi : V^{\flat} \to \RR$ by
\begin{equation*}
\Phi(p,t) = \phi(p)+at
\end{equation*}
Then $\Phi$ is $(\xi, 1)$-compatible.

Let $c \in ( \phi_{max} , \phi_{min}+a)$. The cross section
$V_c^{\flat}$ consists of oriented edges of type $(p,0)(p,1)$
and we can naturally identify it with $V$, the set of vertices of $\Gamma$. 
We can make
$V_c^{\flat}$ into a graph, $\Gamma_c^{\flat}$, by joining the points
corresponding to 
$(p,0)(p,1)$ and $(q,0)(q,1)$ if and only if $p$ and $q$ 
are joined by an edge of $\Gamma$.

\smallskip

\begin{theorem}
\label{th:cohcros=cohske} If $c \in ( \phi_{max} , \phi_{min}+a)$ then
\begin{equation*}
H(\Gamma_c^{\flat}, \alpha_c^{\flat}) \simeq H(\Gamma,\alpha)
\end{equation*}
\end{theorem}

\begin{proof}

There is a natural isomorphism of linear spaces
$\fg^* \to (\fg \oplus \RR)_{(\xi,1)}^* \subset \fg^* \oplus \RR^*$,
given by
$$ \sigma \longrightarrow (\sigma, -\sigma(\xi) \cdot {\bf 1}).$$
This induces an isomorphism of rings
\begin{equation*}
\rho : \SS((\fg\oplus\RR)_{(\xi,1)}^*) \to \SS(\fg^*).
\end{equation*}
and, together with the identification $V_c^{\flat} \simeq V$,
an isomorphism 
\begin{equation}
\label{eq:rhostar}
\rho_* : \mbox{Maps}(V_c^{\flat}, \SS((\fg\oplus\RR)_{(\xi,1)}^*)) 
\to \mbox{Maps}(V,\SS(\fg^*)).
\end{equation}

We will show that \eqref{eq:rhostar} restricts to a bijection
\begin{equation*}
\rho_* : H(\Gamma_c^{\flat}, \alpha_c^{\flat}) \to H(\Gamma,\alpha).
\end{equation*}

First, let $g \in H(\Gamma_c^{\flat}, \alpha_c^{\flat})$.
Let $e$ be an edge of $\Gamma$, with endpoints $p=i(e)$ and $q=t(e)$ 
and let $\Gamma_e$ be the subgraph of $\Gamma^{\flat}$ with vertices
$(p,0), (p,1), (q,0)$ and $(q,1)$. Consider 
the Thom class of $\Gamma_e$ in $\Gamma^{\flat}$. 
This, by definition, is the map
$$h_e : V^{\flat} \to \SS((\fg\oplus\RR)^*)$$ 
that is 0 at vertices not in
$\Gamma_e$ and, at each vertex, $v$, of $\Gamma_e$, 
is equal to the product of the values of $\alpha^{\flat}$ on
edges at $v$ which are not edges of $\Gamma_e$. Then
\begin{equation*}
\frac{\rho_*(g)(p)-\rho_*(g)(q)}{\alpha_{e}} = 
\rho \Bigl( \int_{\Gamma_c^{\flat}} g \K_c(h_e) \Bigr) \in \SS(\fg^*).
\end{equation*}
Since this is true for all edges of $\Gamma$, it follows that 
$\rho_*(g) \in H(\Gamma,\alpha)$.

Conversely, let $g \in H(\Gamma,\alpha)$ and 
$h \in H(\Gamma^{\flat}, \alpha^{\flat})$. The projection on the
first factor $\fg^* \oplus \RR^* \to \fg^*$ induces a map 
$\pi_1 : \SS((\fg \oplus \RR)^*) \to \SS(\fg^*)$. 
Let 
$$h_0=  \pi_1 \circ h |_{\Gamma \times\{ 0 \} };$$
then $h_0 \in H(\Gamma,\alpha)$ and therefore $gh_0 \in H(\Gamma,\alpha)$.

A direct computation shows that
\begin{equation*}
\rho \Bigl( \int_{\Gamma_c^{\flat}} \rho_*^{-1}(g) \K_c(h) \Bigr) =
%\sum_{p \in V} \frac{g(p)h_0(p)}{\prod_i \alpha_{p,e_i}}.
\int_{\Gamma} gh_0 \; .
\end{equation*}

Theorem 2.2 in \cite{GZ1} implies that the right hand side is an element of 
$\SS(\fg^*)$ and hence
\begin{equation*}
\int_{\Gamma_c^{\flat}} \rho_*^{-1}(g) \K_c(h) \in
\SS((\fg \oplus \RR)_{(\xi,1)}^*),
\end{equation*}
for all $h \in H(\Gamma,\alpha)$. Therefore 
$$\rho_*^{-1}(g) \in H(\Gamma_c^{\flat}, \alpha_c^{\flat}),$$ 
which proves the theorem.
\end{proof}

\smallskip

\begin{center}
\section{The changes in cohomology}
\label{sec:6}
\end{center}

\smallskip

Let $c$ and $c'$ be regular values of $\phi$ such that $c < c'$ and such that 
there is exactly one vertex, $p$, with $c < \phi(p) < c'$. Let $r$ be the 
index of $p$, and let $s=d-r$. 
Let $e_i$, $i=1,..,r$ be the edges, $e$, issuing from $p$, for 
which $\alpha_{p,e}(\xi) <0$ and let $e_a$, $a=r+1, .., d$ be 
the other edges issuing from $p$. 
We assume that $1 \leq r \leq d-1$, so that $d-1 \geq s \geq 1$.  
Let $\Delta_c = \{ e_i \; ; \; i=1,..,r \} \subset V_c$ and
$\Delta_{c'} = \{ e_a \; ; \; a=r+1,..,d \} \subset V_{c'}$. Define
\begin{equation*}
V^{\#} = (V_c - \Delta_c) \cup ( \Delta_c \times \Delta_{c'})
\end{equation*}
and a map 
\begin{equation*}
\pi_c : V^{\#} \to V_c
\end{equation*}
which is the identity on $V_c - \Delta_c$ and the projection on the 
first factor on $\Delta_c \times \Delta_{c'}$. Note that 
\begin{equation*}
V^{\#} = (V_{c'} - \Delta_{c'}) \cup ( \Delta_c \times \Delta_{c'})
\end{equation*}
and that there exists a similar map 
\begin{equation*}
\pi_{c'} : V^{\#} \to V_{c'}.
\end{equation*}
We define the pull-back maps
\begin{equation*}
\pi_c^* : \mbox{Maps}(V_c, \SS(\fg_{\xi}^*)) \to 
\mbox{Maps}(V^{\#}, \SS(\fg_{\xi}^*))
\end{equation*}
and
\begin{equation*}
\pi_{c'}^* : \mbox{Maps}(V_{c'}, \SS(\fg_{\xi}^*)) \to 
\mbox{Maps}(V^{\#}, \SS(\fg_{\xi}^*)).
\end{equation*}
The projections of $\Delta_c \times \Delta_{c'}$ onto factors define 
pull-back maps
\begin{equation*}
\mbox{pr}_1^* : \mbox{Maps}(\Delta_c, \SS(\fg_{\xi}^*)) \to 
\mbox{Maps}(V^{\#}, \SS(\fg_{\xi}^*))
\end{equation*}
and
\begin{equation*}
\mbox{pr}_2^* : \mbox{Maps}(\Delta_{c'}, \SS(\fg_{\xi}^*)) \to 
\mbox{Maps}(V^{\#}, \SS(\fg_{\xi}^*)),
\end{equation*}
by extending with 0 outside $\Delta_c \times \Delta_{c'}$. Since
the maps $\pi_c^*, \pi_{c'}^*, \mbox{pr}_1^*$ and $\mbox{pr}_2^*$
are injective, we can regard the sets, 
$\mbox{Maps}(V_c, \SS(\fg_{\xi}^*))$, 
$\mbox{Maps}(V_{c'}, \SS(\fg_{\xi}^*))$,  
$\mbox{Maps}(\Delta_c, \SS(\fg_{\xi}^*))$ and 
$\mbox{Maps}(\Delta_{c'}, \SS(\fg_{\xi}^*))$ as subsets of 
$\mbox{Maps}(V^{\#}, \SS(\fg_{\xi}^*))$.

Let 
\begin{align*}
\alpha_{e_i} & =  m_i(x-\beta_i(y))&  \mbox{and} \\
\alpha_{e_a} & =  m_a(x-\beta_a(y)),&
\end{align*}
with $m_i < 0 < m_a$ and  $\beta_i(y), \beta_a(y) \in \fg_{\xi}^*$,
for $i=1,..,r$ and $a=r+1,..,d$. Consider the maps
\begin{align*}
\tau_c : \; & \Delta_c \to \fg_{\xi}^*,   &\tau_c(e_i)&=\beta_i(y) ;  \\
\tau_{c'} :\;  & \Delta_{c'} \to \fg_{\xi}^* , & \tau_{c'}(e_a) 
&=\beta_a(y) , \mbox{and}\\
\tau^{\#} : \; & \Delta_c \times \Delta_{c'} \to \fg_{\xi}^*  , & 
\tau^{\#}(e_i, e_a) &= \beta_i(y) - \beta_a(y). 
\end{align*}

\begin{remark*}
The fact that $\xi$ is generic implies that $\tau^{\#}(e_i,e_a) \neq 0$.
\end{remark*}

Let $H(\Delta_c, \tau_c)$ be the ring associated to the finite set $\Delta_c$
and the injective function $\tau_c : \Delta_c \to \fg_{\xi}^*$ 
(see the Appendix, Definition \ref{def:abscoh}).

\begin{lemma}
\label{lem:step2}
If $f \in H(\Gamma,\alpha)$ then $f |_{\Delta_c} \in H(\Delta_c, \tau_c)$.
\end{lemma}

\begin{proof}
The proof of this ``localization'' theorem is similar to that of the 
localization theorem in Section \ref{sec:localiz}. 
Let $f_0=f|_{\Delta_c}$ be the restriction of $f$ to $\Delta_c$.
To show that $f_0 \in H(\Delta_c, \tau_c)$ we need to show that
\begin{equation}
\label{eq:intf0}
\int_{\Delta_c} f_0 h(\tau_c) \in \SS(\fg_{\xi}^*)
\end{equation}
for all $h \in \SS(\fg_{\xi}^*)[Y]$, and just as in the proof of 
Theorem \ref{th:new5}, we will show that this can be ``localized'' 
to analogous statements about the integrals of $f_0h(\tau_c)$ over the 
hyperedges of the hypergraph $\Delta_c$. Here are the details.

Let $\MM = \{ \beta_i - \beta_j \; ; \; i \neq j \} \subset \fg_{\xi}^*$.
For $\gamma \in \MM$, let $\MM_{\gamma}$ be the multiplicative system in 
$\SS(\fg_{\xi}^*)$ generated by elements of $\MM$ which are not collinear 
with $\gamma$ and let $\SS(\fg_{\xi}^*)_{\gamma}$ be the corresponding
localized ring. We will show that
\begin{equation}
\label{eq:redgamma}
\int_{\Delta_c} f_0 h(\tau_c) \in \SS(\fg_{\xi}^*)_{\gamma}
\end{equation}
for each $\gamma \in \MM$; this will imply \eqref{eq:intf0}.

Let $\gamma = \beta_i -\beta_j \in \MM$. We join two points $v_k$ and $v_l$ 
of $\Delta_c$ if $\gamma$ and $\beta_k - \beta_l$ are collinear. 

Let $(\Delta_c)_{\gamma}^0 = \{ v_{j_1},..,v_{j_k} \} $ be a non-trivial 
connected component of this new graph. Then $(\Delta_c)_{\gamma}^0$ is a 
complete graph. Let $\fh^*$ be the two-dimensional subspace of $\fg^*$ 
generated by $\alpha_{j_1}$ and $\alpha_{j_2}$ and let 
$\Gamma_{\fh^*}^0$ be the connected component of $\Gamma_{\fh^*}$ 
that contains $p$. Then
\begin{equation*}
(V_{\fh^*}^0)_c \cap \Delta_c = (\Delta_c)_{\gamma}^0.
\end{equation*}

We now use our hypothesis : $\Gamma_{\fh^*}^0$ admits a generating family.
Then the same is true if we change $\xi$ to $-\xi$;
let $\tau_{\fh^*}^{(p)} \in H(\Gamma_{\fh^*},\alpha)$ be a 
generating class at $p$ with respect to $o_{-\xi}$ and define
\begin{equation*}
\label{eq:class}
\psi = \psi_{\fh^*} \tau_{\fh^*}^{(p)} \in H(\Gamma,\alpha),
\end{equation*}
where $\psi_{\fh^*}$ is the Thom class of $\Gamma_{\fh^*}^0$ in $\Gamma$.
Then $\psi$ is supported on the flow-down of $p$ in $\Gamma_{\fh^*}^0$ and
\begin{equation*}
\label{eq:psip}
\psi(p) = \prod_{t \not \in \{j_1,..,j_k \}} \alpha_{e_t}.
\end{equation*}

For every $P \in \SS(\fg^*) \simeq \SS(\fg_{\xi}^*)[Y]$ we have 
\begin{equation}
\label{eq:inton0}
\int_{(\Delta_c)_{\gamma}^0} fP(\tau_c) = - 
(\prod_{t \in \{j_1,...,j_k \}}m_t) 
\int_{\Gamma_c} f\K_c(P\psi);
\end{equation}
since $f \in H(\Gamma_c,\alpha_c)$, \eqref{eq:inton0} implies that
$f|_{(\Delta_c)_{\gamma}^0} \in H((\Delta_c)_{\gamma}^0, \tau_c)$.

Let $R \in \SS(\fg_{\xi}^*)[X]$ be the polynomial
\begin{equation*}
R(X)=\prod_{t \not \in \{ j_1,..,j_k \}} (X-\beta_t)
\end{equation*}
and let $T_0 \in \SS(\fg_{\xi}^*)[X_1,...,X_{k-1}]$ be given by
\begin{equation*}
T_0 = \prod_{l=1}^{k-1} R(X_l).
\end{equation*}
Then $T_0 \in (\SS(\fg_{\xi}^*)[X_1,...,X_{k-1}])^{\Sigma_{k-1}}$, hence,
according to Lemma \ref{lem:2} of the Appendix, there exists 
$T \in (\SS(\fg_{\xi}^*)[X_1,...,X_k])^{\Sigma_k}[Y]$ such that
$T_0=T(X_k)$. In other words, by inserting $X_k$ for $Y$ in this polynomial 
of $k+1$ variables we get back our original polynomial in $k-1$ variables. 
Using this we deduce that
\begin{equation}
\label{eq:intongammac}
\sum_{l=1}^k \frac{f(e_{j_l})h(\beta_{j_l})}
{\prod_{t \neq j_l}(\beta_{j_l}-\beta_t)} = 
\frac{1}{ \prod_{l=1}^k R(\beta_{j_l})} 
\int_{(\Delta_c)_{\gamma}^0} (\Psi h)(\tau_c)f,
\end{equation}
where $\Psi \in \SS(\fg_{\xi}^*)[Y] \simeq \SS(\fg^*)$ is the polynomial 
obtained by replacing $X_l$ with $\beta_{j_l}$, for all $l=1,.,k$.

Since $ \prod_{l=1}^k R(\beta_{j_l}) \in \MM_{\gamma}$, using
\eqref{eq:inton0} we conclude that the left hand side of 
\eqref{eq:intongammac} is an element of 
$\SS(\fg_{\xi}^*)_{\gamma}$. Now \eqref{eq:redgamma} follows from the 
fact that the integral of $f_0h(\tau_c)$ is a sum of terms of this form.
\end{proof}

\begin{lemma}
\label{lem:step3}
For every $f \in H(\Gamma_c, \tau_c)$ and every 
$f_i \in H(\Delta_c, \alpha_c)$, $i=1,..,s-1$, there exist unique 
$f' \in H(\Gamma_{c'}, \alpha_{c'})$ and 
$f_j' \in H(\Delta_{c'}, \tau_{c'})$, $j=1,..,r-1$, such that
\begin{equation}
\label{eq:change}
f' + \sum_{j=1}^{r-1} (\tau^{\#})^jf_j' =
f + \sum_{i=1}^{s-1} (\tau^{\#})^if_i.
\end{equation}
\end{lemma}

\begin{proof}

{\bf Uniqueness}:
It is clear that $f'$ should be equal to $f$ 
on $\Gamma_c - \Delta_c = \Gamma_{c'} - \Delta_{c'}$. If we restrict
\eqref{eq:change} to each set of the form $\Delta_c \times \{ e_a \}$, 
we get a linear system whose matrix is a non-singular Vandermonde matrix; 
therefore, the solution is unique. Hence there exist unique
$f' \in \mbox{Maps}(\Gamma_{c'}, Q(\fg_{\xi}^*))$ and 
$f_j' \in \mbox{Maps}(\Delta_{c'}, Q(\fg_{\xi}^*))$ which satisfy 
\eqref{eq:change}.

{\bf Existence}:
We need to show that $f' \in H(\Gamma_{c'}, \alpha_{c'})$ and 
$f_j' \in H(\Delta_{c'}, \tau_{c'})$, where $f$ and $f_j'$ are the ones 
obtained above.

Using Lemma \ref{lem:invVan} and Corollary \ref{cor:12.1} of the Appendix, 
we deduce that
\begin{equation*}
f_j' = \sum_{m'} \Bigl( \sum_{m,k} P_{m,m',k}(y) 
(\int_{\Delta_c} \tau_c^mf_k) \Bigr) (\tau_{c'})^{m'},
\end{equation*}
where $P_{m,m',k}(y) \in \SS(\fg_{\xi}^*)$;
since $f_k \in H(\Delta_c, \tau_c)$ for all $k=0,..,s-1$, 
we can use Lemma \ref{lem:14.1} to conclude that 
$f_j' \in H(\Delta_{c'}, \tau_{c'})$.

Let $h \in H(\Gamma,\alpha)$. Again, we use 
Lemma \ref{lem:invVan} and Corollary \ref{cor:12.1} to obtain that
\begin{equation*}
\label{eq:f}
\int_{\Gamma_{c'}} f'\K_{c'}(h) = \int_{\Gamma_c} f\K_c(h) +
\sum_{i=0}^{s-1} \int_{\Delta_c} P_i(\tau_c) f_i,
\end{equation*}
where $P_i \in \SS(\fg_{\xi}^*)[X]$. We now use that
$f \in H(\Gamma_c, \alpha_c)$, 
$f_i \in H(\Delta_c, \tau_c)$ for all $i=0,..,s-1$, and 
Lemma \ref{lem:14.1} to conclude that
\begin{equation*}
\label{fprim}
\int_{\Gamma_{c'}} f'\K_{c'}(h) \in \SS(\fg_{\xi}^*)
\end{equation*}
for all $h \in H(\Gamma,\alpha)$, {\it i.e.} that 
$f' \in H(\Gamma_{c'}, \alpha_{c'})$.
\end{proof}

We now use this lemma to determine the change in the cohomology of the 
cross-section as we pass over the vertex $p$.

\begin{corollary}
\label{cor:difdim}
For every $m \geq 0$ we have
\begin{equation}
\label{eq:6.1}
\dim H^{m}(\Gamma_{c'}, \alpha_{c'}) = \dim H^{m}(\Gamma_c,\alpha_c) +
\sum_{k=0}^{s-1} \lambda_{m-k,n-1} - \sum_{k=0}^{r-1} \lambda_{m-k,n-1}.
\end{equation}
\end{corollary}

\begin{proof}
When $0 < \sigma(p) < d$, this follows immediately from 
\eqref{eq:change} and \eqref{eq:14.100}. 
When $\sigma(p)=0$ it follows from \eqref{eq:14.100} and
\begin{equation*}
\label{eq:over0}
H(\Gamma_{c'},\alpha_{c'}) = H(\Gamma_c,\alpha_c) \oplus 
H(\Delta_{c'}, \alpha_{c'})
\end{equation*}
and when $\sigma(p)=d$ it follows from \eqref{eq:14.100} and 
\begin{equation*}
\label{eq:overd}
H(\Gamma_{c},\alpha_{c}) = H(\Gamma_{c'},\alpha_{c'}) \oplus 
H(\Delta_{c}, \alpha_{c}). \qquad \qed
\end{equation*}
\renewcommand{\qed}{}
\end{proof} 

\smallskip

\begin{center}
\section{Adding-up dimensions}
\label{sec:7}
\end{center}

\smallskip

To compute the dimension of $H^m(\Gamma,\alpha)$, we will apply 
Corollary \ref{cor:difdim} several times to the cross-sections of 
$\Gamma^{\flat} = \Gamma \times L$. 

Let $c_0 > \phi_{min}$ such that there is only one vertex, $(p,t)$, of 
$\Gamma^{\flat} =\Gamma \times L$, with $\Phi(p,t) < c_0$. Then
\begin{equation*}
\label{eq:firststep}
\dim H^{m}(\Gamma_{c_0}^{\flat}, \alpha_{c_0}^{\flat}) = 
\sum_{k=0}^{d} \lambda_{m-k,n}.
\end{equation*}
If $c_1 \in ( \phi_{max} , \phi_{min}+a)$ then, according to 
Theorem \ref{th:cohcros=cohske}, we get
\begin{equation*}
\label{laststep}
\dim H^{m}(\Gamma_{c_1}^{\flat}, \alpha_{c_1}^{\flat}) = 
\dim H^m(\Gamma,\alpha).
\end{equation*}

Now let $c_0 < a < b < c_1$ such that there is exactly one vertex of 
$\Gamma$, $p$, such that $a < \Phi(p,0) < b$.  If the index
of $p$ in $\Gamma$ is $\sigma(p)=r$ then the index of $(p,0)$ in 
$\Gamma^{\flat}$ is also $r$. Thus we can apply 
\eqref{eq:6.1} to obtain
\begin{equation*}
\dim H^{m}(\Gamma_b^{\flat}, \alpha_b^{\flat}) = 
\dim H^{m}(\Gamma_a^{\flat}, \alpha_a^{\flat}) +
\sum_{k=1}^{d-r} \lambda_{m-k,n} - \sum_{k=1}^{r-1} \lambda_{m-k,n}.
\end{equation*}
Adding together these changes we get
\begin{equation*}
\label{eq:brack}
\dim H^m(\Gamma,\alpha)  =   
\sum_{p \in V} 
\Bigl( \sum_{k=0}^{d-\sigma(p)}\lambda_{m-k} - 
\sum_{k=0}^{\sigma(p)-1} \lambda_{m-k} \Bigr)
\end{equation*}

The minimum value for $k$ is 0 and the maximum is $d$; $\lambda_{m-k}$
appears in the first sum when $\sigma(p) \leq d-k$ and in the second one when
$\sigma(p) \geq k+1$. Therefore 
\begin{equation*}
\dim H^m(\Gamma,\alpha)  =  
\sum_{k=0}^d \Bigl(
\sum_{l=0}^{d-k} b_{l}(\Gamma) - 
\sum_{l=k+1}^{d} b_{l}(\Gamma) \Bigr) \lambda_{m-k}. 
\end{equation*}
Because $b_{d-l}(\Gamma) = b_{l}(\Gamma)$ (see \eqref{eq:poin}),
the expression in bracket reduces to $b_{k}(\Gamma)$ and therefore
\begin{equation*}
\dim H^m(\Gamma,\alpha) = \sum_{k=0}^d b_{k}(\Gamma) \lambda_{m-k,n} .
\end{equation*}
This concludes the proof of Theorem \ref{th:mainres}.

\begin{remark*}
The results of \cite[Section~2.5.2]{GZ2} are valid in this case
as well. In particular, the dimension of $H^{m}(\Gamma_c, \alpha_c)$ 
is the same 
as the dimension of the image of the Kirwan map 
\begin{equation*}
\K_c : H^m(\Gamma,\alpha) \to \mbox{Maps}(V_c, \SS^m(\fg_{\xi}^*)).
\end{equation*}
Since 
$\mbox{im}(\K_c) \subset H(\Gamma_c, \alpha_c)$, it follows that 
the Kirwan map is surjective. Thus, in particular, 
$H(\Gamma_c, \alpha_c)$ is not only an 
$\SS(\fg_{\xi}^*)$-module, but also a ring.
\end{remark*}

\smallskip

\setcounter{section}{8}

\begin{center}
\section*{Appendix}
\label{sec:4}
\end{center}

\smallskip

\begin{lemma} \label{lem:1}

\begin{equation} \label{eq:12.1} 
\sum_{k=1}^m \frac{X_k^N}{ \prod_{j \neq k}(X_k-X_j)} =
\sum_{\substack{i_1 + .. +i_m=N-m+1 \\ i_1,..,i_m \geq 0}}
X_1^{i_1} \cdots X_m^{i_m}.
\end{equation}
\end{lemma}

\begin{proof}

Consider the decomposition in partial fractions
\begin{equation*}
\frac{Z^N}{(Z-X_1) \cdots (Z-X_m)} = Q(Z) + 
\sum_{k=1}^m \frac{X_k^N}{\prod_{j \neq k}(X_k-X_j)} \frac{1}{Z-X_k},
\end{equation*}
where $Q(Z)$ is a polynomial in $Z$. We use the expansion
\begin{equation}
\label{eq:12.3}
\frac{1}{Z-X_k}= \sum_{i_k=0}^{\infty} X_k^{i_k}Z^{-1-i_k},
\end{equation}
to obtain that:
\begin{equation*}
Z^{N-m}\prod_{k=1}^m\Bigl( \sum_{i_k=0}^{\infty} X_k^{i_k}Z^{-i_k}\Bigr) =
Q(Z) + \frac{1}{Z}
\sum_{k=1}^m \Bigl( \frac{X_k^N}{\prod_{j \neq k}(X_k-X_j)} 
\sum_{l=0}^{\infty} X_k^lZ^{-l}\Bigr).
\end{equation*}
Now the formula \eqref{eq:12.1} follows by comparing coefficients of 
$Z^{-1}$ on both sides.
\end{proof}

\begin{corollary}
\label{cor:12.1}
Let R be a ring and $P \in R[X_1,..,X_m][Y]$. Then
\begin{equation*}
\sum_{k=1}^m \frac{P(X_k)}{ \prod_{j \neq k}(X_k-X_j)} \in R[X_1,..,X_m].
\end{equation*}
\end{corollary}

\begin{lemma}
\label{lem:2}
Let $P_0 \in (\CC[X_1,..,X_{m-1}])^{\Sigma_{m-1}}$ be a symmetric
polynomial. Then there exists 
$P \in (\CC[X_1,..,X_m])^{\Sigma_m}[Y]$ such that $P_0=P(X_m)$.
\end{lemma}

\begin{proof}
It suffices to show that the lemma is true if $P_0$ is a 
fundamental symmetric polynomial in $X_1,..,X_{m-1}$ and for this we can 
use an inductive argument.

We can also start with the identity
\begin{equation*} 
(Z+X_1) \cdots (Z+X_{m-1}) = \frac{(Z+X_1)\cdots (Z+X_m)}{Z} 
\frac{1}{1+X_mZ^{-1}}
\end{equation*}
and we use \eqref{eq:12.3} to deduce that
\begin{equation*} 
(Z+X_1) \cdots (Z+X_{m-1}) = (\sum_{k=0}^m \sigma_k Z^{m-k})
(\sum_{j=0}^{\infty} (-1)^j X_m^j Z^{-j-1}),
\end{equation*}
where $\sigma_k \in (\CC[X_1,..,X_m])^{\Sigma_m}$ is the 
fundamental symmetric polynomial of degree $k$. The lemma follows by 
comparing coefficients of $Z^k$ on both sides.
\end{proof}

\bigskip

\begin{lemma}
\label{lem:invVan}
Let $R=\ZZ[X_1,...,X_m]$ and 
\begin{equation*}
A= 
\begin{pmatrix}
1 & X_1 & \dots & X_1^{m-1} \\
1 & X_2 & \dots & X_2^{m-1} \\
\hdotsfor{3} \\
1 & X_m & \dots & X_m^{m-1} 
\end{pmatrix}
\end{equation*}
be a Vandermonde matrix. Then $A$ is invertible in the space of matrices 
with entries in the quotient field of $R$ and for each $i=1,..,m$, there 
exists $P_i \in R[Y]$ such that 
\begin{equation*}
(A^{-1})_{ij} = \frac{P_i(X_j)}{\prod_{k\neq j}(X_j-X_k)}.
\end{equation*}
Moreover, for $i=1$ we have:
\begin{equation}
\label{eq:13.3}
(A^{-1})_{1j} = \prod_{k\neq j}\Bigl(\frac{-X_k}{X_j-X_k}\Bigr).
\end{equation}
\end{lemma}

\begin{proof}

Since
\begin{equation}
\label{eq:13.4}
\det A = \prod_{1 \leq k < l \leq m}(X_l - X_k) \neq 0,
\end{equation}
the matrix $A$ is invertible; the entries of the inverse are
\begin{equation}
\label{eq:13.5}
(A^{-1})_{ij} = (-1)^{i+j} \frac{\det A_{ji}}{\det A},
\end{equation}
where $A_{ji}$ is obtained from $A$ by removing the $j^{th}$ row and
$i^{th}$ column.

We start with the identity
\begin{equation}
\label{eq:13.6}
(Z-X_1) \cdots (Z-X_{j-1})(Z-X_{j+1}) 
\cdots (Z-X_m) = 
\sum_{k=0}^{m-1} (-1)^{m-1-k} \sigma_{m-1-k}^j Z^k,
\end{equation}
where $\sigma_{m-1-k}^j$ is the fundamental symmetric 
polynomial of degree $m-1-k$ in variables $X_1,..,X_{j-1}, X_{j+1}, .., X_m$.
Since the left hand side of \eqref{eq:13.6} is 0 for $Z=X_l$, $l \neq j$, 
we get:
\begin{equation}
\label{eq:13.7}
 X_l^{m-1} =
\sum_{k=0}^{m-2} (-1)^{m-k} \sigma_{m-k-1}^j X_l^k .
\end{equation}
The left hand side of \eqref{eq:13.7} is the last column of $A_{ji}$; 
using \eqref{eq:13.7} and basic properties of determinants we deduce that
\begin{equation}
\label{eq:13.8} 
\det A_{ji} = \sigma_{m-i}^j 
\prod_{\substack{1 \leq k < l \leq m \\ k \neq j \neq l}}
(X_l-X_k).
\end{equation}
Using \eqref{eq:13.4} and \eqref{eq:13.8} in  \eqref{eq:13.5} we 
obtain that
\begin{equation}
\label{eq:13.9}
(A^{-1})_{ij} = \frac{(-1)^{n+i} \sigma_{n-i}^j}
{\prod_{k\neq j}(X_j -X_k)}.
\end{equation}
We now use Lemma \ref{lem:2} to finish the proof. 
Since $\sigma_{n-1}^j= X_1 \cdots X_{j-1}X_{j+1} \cdots X_m$, 
for $i=1$, \eqref{eq:13.9} becomes \eqref{eq:13.3}.
\end{proof}

\begin{remark*}
The result of Lemma \ref{lem:2} can be made more precise. 
By induction on $k$ we get that
\begin{equation}
\label{eq:skj}
\sigma_{k}^j = \sum_{l=0}^k (-1)^l \sigma_{k-l}X_j^l \; . 
\end{equation}
Then \eqref{eq:13.9} becomes:
\begin{equation}
\label{eq:a-1ij}
(A^{-1})_{ij} = 
\frac{ \sum_{k=0}^{n-i} (-1)^{n+i+k} \sigma_{n-i-k} X_j^k}
{\prod_{k\neq j}(X_j -X_k)}.
\end{equation}
\end{remark*}

\bigskip

Let $W$ be an $n$-dimensional vector space, 
$\Delta = \{ v_1,...,v_d \}$ be a finite set and
$\tau : \Delta  \to W$ be an injective function.
Let $\SS(W)$ be the symmetric algebra and $Q(W)$ be the quotient field of
$\SS(W)$. 

We define an integral operation 
$\int_{\Delta} : \mbox{Maps}(\Delta,Q(W)) \to Q(W)$ by
\begin{equation}
\label{eq:intonfinite}
\int_{\Delta} g = \sum_{k=1}^d \frac{g(v_k)}
{\prod_{j \neq k}(\tau(v_k)-\tau(v_j))}.
\end{equation}
\begin{definition}
\label{def:abscoh}
We define $H(\Delta,\tau)$ to be the space of all maps $g$ satisfying the 
condition 
\begin{equation*}
\int_{\Delta} g P(\tau) \in \SS(W), \; \;  \mbox{ for all } P \in \SS(W)[Y].
\end{equation*}
\end{definition}

\begin{lemma}
\label{lem:14.1}
A map $g : \Delta \to Q(W)$ is in $H(\Delta,\tau)$ if and only if 
there exist $g_0,...,g_{d-1} \in \SS(W)$ such that
\begin{equation}
\label{eq:14.1}
g=\sum_{k=0}^{d-1} g_k\tau^k.
\end{equation}
\end{lemma}

\begin{proof}

The fact that every map of the form \eqref{eq:14.1} is in 
$H(\Delta,\tau)$ is a direct consequence of Corollary \ref{cor:12.1}. 

To show that every element of $H(\Delta,\tau)$ can be written as in 
\eqref{eq:14.1} we proceed as follows: 
We can regard \eqref{eq:14.1} as a linear system whose matrix is a 
Vandermonde matrix; hence there are $g_0,..,g_{d-1} \in Q(W)$
such that \eqref{eq:14.1} is true. Moreover, we can use Lemma 
\ref{lem:invVan} to deduce that
\begin{equation*}
g_i = \int_{\Delta} g P_{i+1}(\tau).
\end{equation*}
Since $g \in H(\Delta,\tau)$,  we conclude that $g_i \in \SS(W)$.
\end{proof}

\noindent {\bf Remarks:}

\begin{enumerate}
\item $H(\Delta,\tau)$ is a graded ring and an $\SS(W)$-module.
\item If $g \in H(\Delta,\tau)$ then $g \in \mbox{Maps}(\Delta,\SS(W))$.
\item If $g \in H(\Delta,\tau)$ then the polynomials $g_0,..,g_{d-1}$ are
unique.
\item Let 
$H^{m}(\Delta,\tau)= H(\Delta,\tau) \cap \mbox{Maps}(\Delta, \SS^m(W))$.
Then 
\begin{equation}
\label{eq:14.100}
\mbox{dim}H^{m}(\Delta,\tau) = \sum_{k=0}^{d-1} \lambda_{m-k,n}.
\end{equation}
\end{enumerate}

\end{document}